\documentclass[12pt,twoside]{amsart}

\usepackage{amssymb}
\usepackage{amsmath}
\usepackage{bbm}
\usepackage{mathrsfs}
\usepackage{float}
\usepackage{amsfonts, latexsym, dsfont}

\newcommand{\br}{ }
\newcommand{\brr}{, }

\makeatletter
\gdef\th@mychange{\normalfont\slshape
   \def\@begintheorem##1##2{\item
        [\hskip\labelsep \theorem@headerfont ##2. ##1  \,--\!--\!--\!--  ]}%
 \def\@opargbegintheorem##1##2##3{%
   \item[\hskip\labelsep \theorem@headerfont ##2. ##1\ {\upshape(}##3{\upshape)}. \,-----  ]}}
\makeatother

\newtheorem{fac}{Fact}[section]
\newtheorem{facs}{}[subsection]

\newtheorem{lem}[fac]{Lemma}

\theoremstyle{definition}
\newtheorem{ttt}[fac]{}
\newtheorem{defi}[fac]{Definition}

\newtheorem{ttts}[facs]{}
\newtheorem{exas}[facs]{Example}

\theoremstyle{remark}
\newtheorem{notas}[facs]{Notation}
\newtheorem{rem}[fac]{Remark}
\newtheorem{rems}[fac]{Remarks}

\newtheorem{remas}[facs]{Remark}
\newtheorem{remass}[facs]{Remarks}

\newcommand{\calO}{\mathscr{O}}
\newcommand{\calS}{\mathscr{S}}

\newcommand{\bbF}{{\mathbbm F}}
\newcommand{\bbQ}{{\mathbbm Q}}
\newcommand{\bbZ}{{\mathbbm Z}}

\newcommand{\pr}{{\rm pr}}

\newcommand{\bP}{{\bf P}}
\newcommand{\rk}{\mathop{\text{\rm rk}}}
\newcommand{\Frob}{\mathop{\rm Frob}}
\newcommand{\Gal}{\mathop{\rm Gal}}
\newcommand{\Div}{\mathop{\rm Div}}
\newcommand{\Pic}{\mathop{\rm Pic}}
\newcommand{\Br}{\mathop{\rm Br}}
\newcommand{\disc}{\mathop{\rm disc}}
\newcommand{\calpp}{{\mathfrak{p}}}
\newcommand{\calqq}{{\mathfrak{q}}}
\newcommand{\et}{\text{\'et}}

\newcounter{abc}
\newenvironment{abc}{\begin{list}{\rm \alph{abc}) }%
{\usecounter{abc} \leftmargin=0.0pt \labelsep=0.0pt %
\listparindent=0.0pt \labelwidth=0.0pt \parsep=\smallskipamount %
\itemsep=0.0pt \topsep=0.0pt \partopsep=\smallskipamount}}{\end{list}}
\newcounter{iii}
\newenvironment{iii}{\begin{list}{\rm \roman{iii}) }%
{\usecounter{iii} \leftmargin=0.0pt \labelsep=0.0pt %
\listparindent=0.0pt \labelwidth=0.0pt \parsep=\smallskipamount%
 \itemsep=0.0pt \topsep=0.0pt \partopsep=\smallskipamount}}{\end{list}}

\def\rightend#1#2{{%
 \leavevmode\nobreak\hskip .5em plus 1fil
 \penalty600 \hskip 0pt plus -1filll
 \vadjust{}\nobreak\hskip 0pt plus 1filll%
 #1\parfillskip=#2\relax \par}}

\def\eop{\ifmmode\rule[-22pt]{0pt}{1pt}\ifinner\tag*{$\square$}\else\eqno{\square}\fi\else\rightend{$\square$}{0pt}\fi}

\author[Andreas-Stephan Elsenhans]{ by \\ \\ Andreas-Stephan Elsenhans${}^*$ (Bayreuth) \\}
\author[J\"org Jahnel]{\\ J\"org Jahnel (Siegen)}

\date{}

\title{On the computation of the Picard group
for~$K3$~surfaces${}^\ddagger$}


\address{Fachbereich 6 Mathematik \\
Universit\"at Siegen \\
Walter-Flex-Stra\ss e 3 \\
D-57068 Siegen \\
Germany}

\address{Mathematisches Institut \\
Universit\"at Bayreuth \\
D-95440 Bayreuth \\
Germany}

\begin{document}

\maketitle

\section{Introduction}

\renewcommand{\thefootnote}{\fnsymbol{footnote}}
\footnotetext[1]{The first author was partially supported by the Deutsche Forschungsgemeinschaft (DFG) through a funded research~project.}
\footnotetext[3]{The computer part of this work was executed on the Sun Fire V20z Servers of the Gau\ss\ Laboratory for Scientific Computing at the G\"ottingen Mathematisches~Institut. Both~authors are grateful to Prof.~Y.~Tschinkel for the permission to use these machines as well as to the system administrators for their~support.}

\begin{ttt}
In~this note, we will present a method to construct examples of
$K3$~surfaces
of geometric Picard
rank~$1$.
Our~approach is a refinement of that of R.~van~Luijk~\cite{vL}. It~is based on an analysis of the Galois~module structure on \'etale~cohomology. This~allows to abandon the original limitation to cases of Picard
rank~$2$
after reduction
modulo~$p$.
Furthermore,~the use of Galois~data enables us to construct examples which require significantly less computation~time.
\end{ttt}

\begin{ttt}
The Picard group of a $K3$~surface $S$ is a highly interesting~invariant.
In~general, it is
isomorphic to $\bbZ^n$ for some $n=1,\ldots,20$.
The first explicit examples of $K3$~surfaces over $\bbQ$ with geometric
Picard~rank $1$ were constructed by R.~van Luijk~\cite{vL}.
His~method is based on reduction modulo~$p$.
It~works as~follows.\smallskip

\begin{iii}
\item
At a
place~$p$
of good reduction, the Picard
group~$\Pic(S_{\overline{\bbQ}})$
of the surface injects into the
Picard
group~$\smash{\Pic(S_{\overline{\bbF}_{\!p}})}$
of its reduction
modulo~$p$.
\item
On~its part,
$\smash{\Pic(S_{\overline{\bbF}_{\!p}})}$
injects into the second \'etale cohomology group
$\smash{H_\et^2(S_{\overline{\bbF}_{\!p}}\!, \bbQ_l(1))}$.\vspace{0.1mm}
\item
Only roots of unity can arise as eigenvalues of the
Frobenius on the image
of~$\Pic(S_{\overline{\bbF}_{\!p}})$
in~$\smash{H_\et^2(S_{\overline{\bbF}_{\!p}}\!, \bbQ_l(1))}$.
The~number of eigenvalues of this form is therefore an upper bound for the Picard rank
of~$S_{\overline{\bbF}_{\!p}}$.
We~can compute the eigenvalues
of~$\Frob$
by counting the points
on~$S$,
defined
over~$\bbF_{\!p}$
and some finite~extensions.
\end{iii}\smallskip

\noindent
Doing~this for one prime, one obtains an upper bound
for~$\rk \Pic(S_{\overline{\bbF}_{\!p}})$
which is always~even. The~Tate conjecture asserts that this bound is actually~sharp.

When~one wants to prove that the Picard rank
over~$\overline{\bbQ}$
is, in fact, equal
to~$1$,
the best which could happen is to find a prime that yields an upper bound
of~$2$.
There~is not much hope to do better when working with a single~prime.\smallskip

\begin{iii}
\item[iv) ]
In~this case, the assumption that the surface would have
Picard~rank~$2$
over~$\overline{\bbQ}$
implies that the discriminants of both Picard groups,
$\Pic(S_{\overline{\bbQ}})$
and~$\Pic(S_{\overline{\bbF}_{\!p}})$,
are in the same square~class. Note~here that reduction
modulo~$p$
respects the intersection~product.
\item[v) ]
When one combines information from two primes, it may happen
that we get the rank
bound~$2$
at both places but different
square classes for the discriminant do~arise. Then,~these data are incompatible
with Picard
rank~$2$
over~$\overline{\bbQ}$.

On~the other hand, there is a non-trivial divisor known explicitly.
Altogether,
rank~$1$
is~proven.
\end{iii}
\end{ttt}

\begin{rem}
This~method has been applied by several authors in order to construct
$K3$~surfaces with
prescribed Picard rank~\cite{vL,Kl,EJ1}.
\end{rem}

\begin{ttt}
{\bf The refinement.}
In this note, we will refine van~Luijk's~method. Our~idea is the~following. We~do not look at the ranks,~only. We~analyze the Galois module structures on the Picard~groups,~too. The~point here is that a Galois~module typically has submodules by far not of every~rank.

As~an example, we will construct
$K3$~surfaces
of geometric
Picard~rank~$1$
such that the reduction
modulo~$3$
has geometric Picard rank~$4$
and the reduction
modulo~$5$
has geometric Picard 
rank~$14$.
\end{ttt}

\begin{rem}
This~work continues our investigations on Galois module structures on the Picard~group. In~\cite{EJdisc,EJ2,EJ3}, we constructed
cubic~surfaces~$S$
over~$\bbQ$
with prescribed Galois module structure
on~$\Pic(S)$.
\end{rem}

\section{The Picard group as a Galois~module}
\label{strat}

\begin{ttt}
Let
$K$
be a field and
$S$
an algebraic surface defined
over~$K$.
Denote~by
$S$
the
$\bbQ$-vector
space~$\Pic(S_{\overline{K}}) \!\otimes_\bbZ\! \bbQ$.
On~$S$,
there is a natural
\mbox{$\Gal(\overline{K}/K)$-operation.}
The~kernel of this representation is a normal subgroup of finite~index.
It~corresponds to a finite Galois extension
$L$
of~$K$.
In~fact, we have a
\mbox{$\Gal(L/K)$-representation}.
%We~will call
%$L$
%the field of definition of the Picard~group.

The~group
$\Gal(\overline{K}/L)$
acts trivially
on~$\Pic(S_{\overline{K}})$.
I.e.,
$$\Pic(S_{\overline{K}}) = \Pic(S_{\overline{K}})^{\Gal(\overline{K}/L)} \, .$$
Within~this,
$\Pic(S_L)$
is, in general, a subgroup of finite~index. Equality~is true under the hypothesis
that~$S(L) \neq \emptyset$.
\end{ttt}

\begin{ttt}
Now suppose
$K$
is a number field and
$\calpp$
is a prime ideal
of~$K$.
We~will denote the residue class field
by~$k$.
Further,~let
$S$
be a
$K3$~surface
over~$K$
with good reduction
at~$\calpp$.
%It is well known that reduction modulo $\calpp$ respects the
%intersection product.
There~is an injection
of~$\Pic(S_{\overline{K}})$
into~$\Pic(S_{\overline{k}})$.
Taking~the tensor product, this yields an inclusion of vector spaces
$\Pic(S_{\overline{K}}) \!\otimes_\bbZ\! \bbQ \hookrightarrow \Pic(S_{\overline{k}}) \!\otimes_\bbZ\! \bbQ$.

Both~spaces are equipped with a Galois~operation.
On~$\Pic(S_{\overline{K}}) \!\otimes_\bbZ\! \bbQ$,
we have
a~$\Gal(L / K)$-action.
On~$\Pic(S_{\overline{k}}) \!\otimes_\bbZ\! \bbQ$,
only
$\smash{\Gal(\overline{k}/k) = \overline{\langle\Frob\rangle}}$~operates.
\end{ttt}

\begin{lem}
The~field extension\/
$L/K$
is unramified
at\/~$\calpp$.\smallskip

\noindent
Proof.
{\em
Let~$D \in \Div(S_L)$
be an arbitrary~divisor. By~good reduction,
$D$
extends to a divisor on a smooth model
$\calS$
over the integer
ring~$\calO_L$.
In~particular, we have the reduction
$D_\calqq$
of~$D$
on the special fiber
$S_\calqq$.
Here,~$\calqq$
is any prime, lying
above~$\calpp$.

$\calO_L/\calqq$~is
a finite extension
of~$k$.
Correspondingly,~there is an unramified
extension~$L^\prime \subset L$
of~$K$.
Good~reduction implies that every divisor
on~$S_\calqq$
lifts
to~$S_{L^\prime}$.
Consequently,~we have a divisor
$D^\prime \in \Div(S_{L^\prime})$
which has the same reduction
as~$D$.

As~intersection products are respected by reduction, we see that the intersection number
of~$D^\prime_L - D$
with any divisor is~zero. The~standard argument from~\cite[Proposition~VIII.3.6.i)]{BPV} implies
$D^\prime_L - D = 0$.
In~other words,
$D$
is defined over an unramified~extension.
}
\eop
\end{lem}

\begin{ttt}
There~is a Frobenius lift
to~$L$
which is unique up to~conjugation. When~we choose a particular
prime~$\calqq$,
lying
above~$\calpp$,
we fix a concrete Frobenius~lift.
Then,~$\Pic(S_{\overline{K}}) \!\otimes_\bbZ\! \bbQ$
becomes a
$\Gal(\overline{k} / k)$-submodule
of~$\Pic(S_{\overline{k}}) \!\otimes_\bbZ\! \bbQ$.
\end{ttt}

\begin{ttt}
{\bf Computability of the Galois representation.}
The~simplest way to understand the
$\Gal(\overline{k} / k)$-representation
on ~$\Pic(S_{\overline{k}}) \!\otimes_\bbZ\! \bbQ$
is to use \'etale~cohomology. Counting~the numbers of points,
$S$~has
over finite extensions
of~$k$,
we can compute the characteristic
polynomial~$\Phi$
of the Frobenius
on~$H_\et^2(S_{\overline{k}}\!, \bbQ_l(1))$.
This~is actually a polynomial with rational, even integer,~coefficients and independent of the choice
of~$l\neq p$~\cite[Th\'eor\`eme~1.6]{De}.

Denote
by~$V_{\text{Tate}}$
the largest subspace
of~$H_\et^2(S_{\overline{k}}\!, \bbQ_l(1))$
on which all eigenvalues of the Frobenius are roots of~unity.
On~the other hand,
let~$P_{\text{conj}}$
be the subgroup
of~$\Pic(S_{\overline{K}})$
generated by the conjugates of all the divisors we know~explicitly.

Then,~we have the following chain of
$\Gal(\overline{k}/k)$-modules,
$$
H_\et^2(S_{\overline{k}}, \bbQ_l(1))
\supseteq V_{\text{Tate}} \supseteq
\Pic(S_{\overline{k}}) \!\otimes_\bbZ\! \bbQ_l \supseteq
\Pic(S_{\overline{K}}) \!\otimes_\bbZ\! \bbQ_l \supseteq
P_{\text{conj}}  \!\otimes_\bbZ\! \bbQ_l \, .
$$
In~an optimal situation, the quotient space
$V_{\text{Tate}} / (P_{\text{conj}} \!\otimes_\bbZ\! \bbQ_l)$
has only finitely many
$\Gal(\overline{k}/k)$-submodules.
This~finiteness condition generalizes the codimension one condition, applied in van Luijk's method, step~v).

Our main strategy will then be as~follows. We~inspect the \mbox{$\Gal(\overline{k}/k)$-submodules}
of~$V_{\text{Tate}} / (P_{\text{conj}} \!\otimes_\bbZ\! \bbQ_l)$.
For~all these, except for the null space, we aim to exclude the possibility that it coincides with
$(\Pic(S_{\overline{K}}) \!\otimes_\bbZ\! \bbQ_l) / (P_{\text{conj}} \!\otimes_\bbZ\! \bbQ_l)$.
\end{ttt}

\begin{rems}
\begin{abc}\item
A~sufficient criterion for
a~$\Gal(\overline{k}/k)$-module
to have only finitely many submodules is that the characteristic polynomial
of~$\Frob$
has only simple~roots. This~fact, although very standard, is central to our~method.
\item
Only~submodules of the
form~$P \!\otimes_\bbZ\! \bbQ_l$
for
$P$
a
\mbox{$\Gal(\overline{k}/k)$-submodule}
of~$\Pic(S_{\overline{k}})$
are possible candidates
for~$\Pic(S_{\overline{K}}) \!\otimes_\bbZ\! \bbQ_l$.
Such~submodules automatically lead to factors
of~$\Phi$
with coefficients
in~$\bbQ$.
\end{abc}
\end{rems}

\begin{defi}
We will call a
$\Gal(\overline{k}/k)$-submodule
of~$H_\et^2(S_{\overline{k}}, \bbQ_l(1))$
{\em admissible\/} if it is a
$\bbQ_l$-sub\-vector~space
and the characteristic polynomial
of~$\Frob$
has rational~coefficients.
\end{defi}

\begin{rem}
In some sense, we apply the van der Waerden criterion to the
representations
of~$\Gal(\overline{\bbQ}/\bbQ)$
on the Picard group and \'etale~cohomology.
\end{rem}

\begin{rem}
In~the practical computations presented below, we will work
with~$\bbQ_l$
instead
of~$\bbQ_l(1)$.
This~is the canonical choice from the point of view of counting points but not for the image of the Picard~group. The~relevant zeroes of the characteristic polynomial of the Frobenius are then those of the form
$q$
times a root of~unity.
\end{rem}

\section{An example}

\subsection{Formulation}

\begin{exas}
\label{examp}
Let
$S \colon w^2 = f_6(x,y,z)$
be a
$K3$~surface
of
degree~$2$
over~$\bbQ$.
Assume~the congruences
$$
f_6 \equiv y^6 + x^4  y^2 - 2 x^2  y^4 + 2 x^5 z + 3 x z^5 + z^6 \pmod 5
$$
and
\begin{eqnarray*}
f_6 &\equiv&
2 x^6 + x^4 y^2 + 2 x^3 y^2 z + x^2 y^2 z^2 + x^2 y z^3 + 2 x^2 z^4 \\
&&{}+ x y^4 z + x y^3 z^2 + x y^2 z^3 + 2 x z^5 + 2 y^6 + y^4 z^2 + y^3 z^3
\pmod 3\, .
\end{eqnarray*}
Then,~$S$
has geometric
Picard~rank~$1$.
\end{exas}

\subsection{Explicit divisors}

\begin{notas}
We will write
$\pr\colon S \to \bP^2$
for the canonical projection.
On~$S$,
there is the ample divisor
$H := \pi^* L$
for
$L$
a line
on~$\bP^2$.
\end{notas}

\begin{ttts}
Let~$C$
be any irreducible~divisor
on~$S$.
Then,~$D := \pi_* C$
is a curve
in~$\bP^2$.
We~denote its degree
by~$d$.
The~projection from
$C$
to~$D$
is generically 2:1 or~1:1. In~the case it is 2:1, we have
$C = \pi^* D \sim dH$.

Thus, to generate a Picard group of
rank~$>\!\! 1$,
divisors are needed which are generically 1:1 over their~projections. This~means,
$\pi^* D$
must be reducible into two components which we call the {\em splittings\/}
of~$D$.

A~divisor
$D$
has a split pull-back if and only if
$f_6$
is a perfect square
on (the normalization
of)~$D$.
A~necessary condition is that the intersection of
$D$
with the ramification locus is a
$0$-cycle
divisible
by~$2$.
\end{ttts}

\subsection{The Artin-Tate conjecture}

\begin{ttts}
The Picard group of a projective variety is equipped with a
$\bbZ$-valued
bilinear form, the intersection~form.
Therefore,~associated
to~$\Pic(S_{\overline{k}})$,
we have its~discriminant, an~integer. The~same applies to every subgroup
of~$\Pic(S_{\overline{k}})$.

For~a
$\bbQ_l$-vector~space 
contained
in~$\Pic(S_{\overline{k}}) \!\otimes_\bbZ\! \bbQ_l$,
the discriminant is determined only up to a factor being a square
in~$\bbQ_l$.
However,~every non-square
in~$\bbQ$
is a non-square
in~$\bbQ_l$
for some suitable
prime~$l \gg 0$.
\end{ttts}

\begin{ttts}
Let us recall the Artin-Tate conjecture in the special case of a
$K3$~surface.\medskip

\noindent
{\bf Conjecture}
(Artin-Tate)
{\bf .}
{\em
Let\/
$Y$
be a\/
$K3$~surface
over\/~$\bbF_{\!q}$.
Denote~by\/
$\rho$
the rank
and by\/
$\Delta$
the discriminant of the Picard group
of\/~$Y$,
defined
over\/~$\bbF_{\!q}$.
Then,
\begin{eqnarray}
\label{AT}
\lim_{T \rightarrow q} \frac{\Phi(T)}{(T-q)^{\rho}} =
q^{21 -\rho} \#\!\Br(Y) |\Delta| \, .
\end{eqnarray}
Here,
$\Phi$
is the characteristic polynomial of the Frobenius
on\/~$H_\et^2(Y_{\overline{\bbF}_{\!q}}\!, \bbQ_l)$.
Finally,~$\Br(Y)$
denotes the Brauer group
of\/~$Y$.
}
\end{ttts}

\begin{remass}
\begin{abc}
\item
The Artin-Tate conjecture allows to compute the square class
of the discriminant of the Picard group over a finite field without any knowledge of explicit~generators.
\item
Observe~that
$\#\!\Br(Y)$
is always a perfect square~\cite{LLR}.
\item
The~Artin-Tate conjecture is proven for most
$K3$~surfaces.
%We~will use it to get good~candidates. Following~Kloosterman~\cite{Kl},
Most~notably, the Tate conjecture implies the Artin-Tate conjecture~\cite{Mi2}. We~will use the Artin-Tate conjecture only in situations where the Tate~conjecture is~true. Thus,~our final result will not depend on
unproven~statements.
\end{abc}
\end{remass}

\subsection{The modulo 3 information}

\begin{ttts}
The sextic curve given by
``$f_6=0$''
has three conjugate conics, each tangent in six~points. Indeed,~note that, for
$$f_3 := x^3 + 2 x^2 y + x^2 z + 2 x y^2 + x y z + x z^2 + y^3 + y^2 z + 2 y z^2 + 2 z^3,$$
the term
$f_6 - f_3^2$
factors into three quadratic forms
over~$\bbF_{\!27}$.
Consequently,~we have three divisors
on~$\bP^2_{\bbF_{\!27}}$
the pull-backs of which~split.
\end{ttts}

\begin{ttts}
Counting~the points
on~$S$
over~$\bbF_{\!3^n}$
for~$n=1, \ldots, 11$
yields the numbers
$-2, -8, 28, 100, 388, 2\,458, 964, -692$, $26\,650$,
$-20\,528$,
and
$-464\,444$
as the traces of the iterated Frobenius
on~$H_\et^2(S_{\overline{\bbF}_3}\!,\bbQ_l)$.
Taking~into account the fact that
$p$~is
a root of the characteristic
polynomial~$\Phi$,
these data determine this polynomial~uniquely,
\begin{eqnarray*}
\Phi(t) & = &t^{22} + 2 t^{21} + 6 t^{20} - 27 t^{18} - 162 t^{17} - 729 t^{16}
- 1458 t^{15} - 2187 t^{14} \\
&&{} + 19683 t^{12}+ 118098 t^{11} + 177147 t^{10} - 1594323 t^8 - 9565938 t^7 \\
&&{} \hspace{0.6cm} - 43046721 t^6- 86093442 t^5 -  129140163 t^4 + 2324522934 t^2 \\
&&{} \hspace{5.8cm} + 6973568802 t \!+\! 31381059609 \, .
\end{eqnarray*}
The~functional equation holds with the plus~sign. We~factorize and~get
\begin{eqnarray*}
\Phi(t) & = &(t - 3)^2
(t^2 + 3 t + 9) \\
&&(t^{18} + 5 t^{17} + 21 t^{16} + 90 t^{15} + 297 t^{14} + 891 t^{13}
+ 2673 t^{12}  + 7290 t^{11} \\
&&{} \hspace{1.3cm} + 19683 t^{10} + 59049 t^9 + 177147 t^8 + 590490 t^7 + 1948617 t^6 \\
&&{} \hspace{1.8cm} + 5845851 t^5  + 17537553 t^4 + 47829690 t^3 + 100442349 t^2 \\
&&{} \hspace{6.8cm} + 215233605 t + 387420489) \, .
\end{eqnarray*}
\end{ttts}

\begin{ttts}
From this, we derive an upper bound
of~$4$
for the rank of the Picard~group.
%But~we can say~more.
%$\Pic(S_{\overline{\bbQ}})$
%leads to a sub-vector~space
%of~$H_\et^2(S_{\overline{\bbQ}}\!, \bbQ_l)$.
%But~this is not only a
%$\bbQ_l$-vector~space.
%It~is
%a~$\smash{\Gal(\overline{\bbQ}/\bbQ)}$-submodule.
%
%When we reduce
%modulo~$p$,
%we pass to a subgroup of the absolute Galois~group.
%Therefore,~the image
%of~$\Pic(S_{\overline{\bbQ}}) \!\otimes_\bbZ\! \bbQ_l$
%is a
%$\Gal(\overline{\bbF}_{\!q}/\bbF_{\!q})$-submodule
%of
%$H_\et^2(S_{\overline{\bbF}_{\!p}}\!, \bbQ_l) \cong H_\et^2(S_{\overline{\bbQ}}\!, \bbQ_l)$.
In~the notation of section~\ref{strat},
$V_{\text{Tate}}$
is 
a~$\bbQ_l$-vector~space
of dimension~four. On~the other hand,
$P_{\text{conj}}$
is generated
by~$H$.
As~$H$
corresponds to one of the factors
$(t-3)$,
the characteristic polynomial of the Frobenius
on~$V_{\text{Tate}} / (P_{\text{conj}} \!\otimes_\bbZ\! \bbQ_l)$
is~$(t - 3)(t^2 + 3 t + 9)$.
It~has only simple~roots.

Consequently,~for each of the dimensions
$1$,
$2$,
$3$,
and~$4$,
there is precisely one admissible
$\Gal(\overline{\bbF}_{\!3}/\bbF_{\!3})$-submodule
in~$H_\et^2(S_{\overline{\bbF}_{\!3}}\!, \bbQ_l)$
containing the Chern class
of~$H$.
%These~are the possible candidates for the image
%of~$\Pic(S_{\overline{\bbQ}}) \!\otimes_\bbZ\! \bbQ_l$.
\end{ttts}

\begin{ttts}
Let us compute the corresponding~discriminants.

\begin{iii}
\item
In the one-dimensional case, we have
discriminant~$2$.
\item
It~turns out that splitting the divisors given by the conics which are six times tangent yields a rank~three
submodule~$M$
of~$\Pic(S_{\overline{\bbF}_{\!3}})$.
Its~discriminant~is
$$\disc M = \det
\left(
\begin{array}{rrr}
-2 & 6 & 0 \\
 6 &-2 & 4 \\
 0 & 4 &-2 \\
\end{array}
\right)
= 96 \, .
$$
Hence,~in the three-dimensional case, the discriminant is in the square class 
of~$6$.
\item
For~the case of dimension four, we may suppose that
$\Pic(S_{\overline{\bbF}_{\!3}})$
is of rank~four.
As~$\smash{\Gal(\overline{\bbF}_3/\bbF_{\!27})}$
acts trivially on 
$\Pic(S_{\overline{\bbF}_{\!3}})$,
the group
$\Pic(S_{\bbF_{\!27}})$
is of rank four,~already. This~means, the Tate conjecture is true
for~$S_{\bbF_{\!27}}$.

We~may compute the square class of the corresponding discriminant according to the Artin-Tate~conjecture. The~result
is~$(-163)$.
\item
Finally,~consider the two-dimensional~case. We~may suppose that
$\Pic(S_{\overline{\bbF}_{\!3}})$
has a
$\smash{\Gal(\overline{\bbF}_{\!3}/\bbF_{\!3}})$-submodule~$N$
which is of rank two and
contains~$H$.
The~corresponding characteristic polynomial of the Frobenius is, necessarily, equal
to~$(t-3)^2$.

On~the other hand, there is the rank~three
submodule~$M$
generated by the splittings of the conics which are six times~tangent. As~the corresponding factors
are~$(t-3)(t^2+3t+9)$,
the modules
$N$
and~$M$
together generate
rank~$4$.
The~Tate conjecture is~true
for~$S_{\bbF_{\!27}}$.
Consequently,~it is true
for~$S_{\bbF_{\!3}}$,~too.

Using~the Artin-Tate conjecture, we can compute the square class of the~discriminant.
It~turns out to
be~$(-489)$.
\end{iii}
\end{ttts}

\begin{remas}
The~Tate conjecture predicts Picard
rank~$2$
for~$S_{\bbF_{\!3}}$.
Let~$C$
be an irreducible divisor linearly independent
of~$H$.
Then,~$C$
is a splitting of a
curve~$D$
of
degree~$d \geq 23$.
Indeed,~$H$
is a genus
$2$
curve.
Hence,~$H^2 = 2$.
For~the discriminant, we find
$-489 \geq 2C^2 - d^2$.
As~$C^2 \geq -2$,
the assertion~follows.

Further,~$D$
is highly singular on the ramification~locus. In~fact, we have
$\smash{C^2 \leq \frac{d^2 - 489}2}$
and~$D^2 = d^2$.
Hence,~going from
$D$
to~$C$
lowers the arithmetic genus by at
least~$\frac{d^2 + 489}4$.
\end{remas}

\subsection{The modulo 5 information}

\begin{ttts}
The sextic curve given~by
``$f_6 = 0$''
has six tritangent~lines. These~are given by
$L_a \colon t \mapsto [1 : t : a]$
where
$a$
is a zero
of~$a^6 + 3 a^5 + 2 a$.
The pull-back of each of these lines splits on the
$K3$~surface
$S_{\overline{\bbF}_{\!5}}$.
\end{ttts}

\begin{ttts}
On~the other hand, counting points yields the following traces of the iterated
Frobenius on~$H_\et^2(S_{\overline{\bbF}_{\!5}}\!, \bbQ_l)$,
$$15, 95, -75, 2\,075, -1\,250, -14\,875, 523\,125, 741\,875, 853\,125, 11\,293\,750 \, .$$
This~leads to the characteristic polynomial
\begin{eqnarray*}
\Phi(t) & = & t^{22} - 15t^{21} + 65t^{20} + 175t^{19} - 3000t^{18} + 11437t^{17} + 10630t^{16} \\
&& {} - 385950t^{15} + 2445250t^{14} - 4530625t^{13} - 38478125t^{12} \\
&& {} + 305656250t^{11} - 566328125t^{10} - 4809765625t^9 + 38207031250t^8 \\
&& {} - 101308593750t^7 - 143457031250t^6 + 2792236328125t^5 \\
&& {} - 14189453125000t^4 + 16400146484375t^3 + 247955322265625t^2 \\
&& {} - 1430511474609375t + 2384185791015625 \\
& = & (t - 5)^2
(t^4 + 5 t^3 + 25 t^2 + 125 t + 625) \\
&& (t^8 - 5 t^7 + 125 t^5 - 625 t^4 + 3\,125 t^3 - 78\,125 t + 390\,625) \\
&& (t^8 - 5 t^7 - 10 t^6 + 75 t^5 - 125 t^4 + 1\,875 t^3 - 6\,250 t^2 - \ldots \\
&& {} \hspace{7.7cm} \ldots - 78\,125 t + 390\,625) \, .
\end{eqnarray*}
Observe here the first two factors correspond to the part of the
Picard group generated by the splittings of the six tritangent~lines.
They~could have been computed directly from the intersection matrix
of these~divisors.
\end{ttts}

\begin{remas}
The knowledge of these two factors allows to compute the characteristic
polynomial only from the numbers of points
over~$\bbF_{\!5},\ldots,\bbF_{\!5^8}$.
Counting~them takes approximately five~minutes when one uses the method described in~\cite[Algorithm~15]{EJ1}.
\end{remas}

\begin{ttts}
Here,~$V_{\text{Tate}}$
is a vector space of
dimension~$14$.
Again,~$P_{\text{conj}}$
is generated
by~$H$.
The~characteristic polynomial of the Frobenius
on~$V_{\text{Tate}} / (P_{\text{conj}} \!\otimes_\bbZ\! \bbQ_l)$~is
\begin{eqnarray*}
(t - 5)(t^4 + 5 t^3 + 25 t^2 + 125 t + 625) \hspace{6.8cm} \\
\hspace{3.5cm} (t^8 - 5 t^7 + 125 t^5 - 625 t^4 + 3\,125 t^3 - 78\,125 t + 390\,625)
\end{eqnarray*}
having only simple~roots. This~shows that, in each of the dimensions
$1$,
$2$,
$5$,
$6$,
$9$,
$10$,
$13$,
and~$14$,
there is precisely one admissible
$\Gal(\overline{\bbF}_{\!5}/\bbF_{\!5})$-submodule
in~$H_\et^2(S_{\overline{\bbF}_{\!5}}\!, \bbQ_l)$
containing the Chern class
of~$H$.
\end{ttts}

\begin{ttts}
For~the cases of low rank, let us compute the square classes of the~discriminant.

\begin{iii}
\item
In the one-dimensional case, we have
discriminant~$2$.
\item
For the two-dimensional case, recall that we know six tritangent lines of the ramification~locus. One~of them,
$L_0$,
is defined
over~$\bbF_{\!5}$.
Splitting~$\pi^* L_0$
yields rank two~alone. For~the discriminant, we find
$$
\det
\left(
\begin{array}{rr}
-2 &  3 \\
 3 & -2 \\
\end{array}
\right)
= -5.
$$
\end{iii}
\end{ttts}

\begin{remas}
Using~the Artin-Tate conjecture, we may compute conditional values for the square classes of the discriminant for the
$6$-
and
$14$-dimen\-sional~modules.
Both~are actually equal
to~$(-1)$.
\end{remas}

\subsection{The situation
over~$\overline{\bbQ}$}

\begin{ttts}
{\em Proof of~\ref{examp}.}
Now we can put everything together and show that the
$K3$~surfaces
described in
Example~\ref{examp}
indeed have geometric Picard~rank~1.

The Picard
group~$\Pic(S_{\overline{\bbQ}}) \!\otimes_\bbZ\! \bbQ_l$
injects as a Galois~submodule into the second \'etale cohomology
groups~$H_\et^2(S_{\overline{\bbF}_{\!5}}\!, \bbQ_l)$
for
$p=3$
and~$5$.
The~modulo~$3$
data show that this module has
\mbox{$\bbQ_l$-dimension}
$1,2,3$
or~$4$.
The~reduction
modulo~$5$
allows the
\mbox{$\bbQ_l$-dimensions}
$1,2,5,6,9,10,13$
and~$14$.
Consequently,~the Picard~rank is either
$1$
or~$2$.

To exclude the possibility of
rank~$2$,
we compare the~discriminants. The~reduction
modulo~$3$
enforces
discriminant~$(-489)$
while the reduction
modulo~$5$
yields
discriminant~$(-5)$.
This~is a~contradiction, e.g.,
for~$l = 17$.
\eop
\end{ttts}

\end{document}